\documentclass[12pt]{article} \textheight=9truein
\newtheorem{lem}{Lemma}
\newtheorem{thm}[lem]{Theorem}
\newtheorem{cor}[lem]{Corollary}
\def\qed{\hspace*{\fill}\vrule height6pt width4pt depth0pt\bigskip}

\long\def\longdelete#1{}

\def\ip{{\rm ip}}
\usepackage{amssymb}

\begin{document}

\title{\bf Isometric path numbers of graphs\thanks{%
                  Supported in part by the National Science Council
                  under grant NSC92-2115-M-002-015.}}
\author{Jun-Jie Pan \\
        Department of Applied Mathematics \\
        National Chiao Tung University \\
        Hsinchu 300, Taiwan \\ \\
        \and
        Gerard J. Chang \\
        Department of Mathematics \\
        National Taiwan University \\
        Taipei 106, Taiwan \\
        E-mail: gjchang@math.ntu.edu.tw
       }
\date{October 3, 2003} 

\maketitle

\begin{abstract}
An isometric path between two vertices in a graph $G$ is a
shortest path joining them. The isometric path number of $G$,
denoted by $\ip(G)$, is the minimum number of isometric paths
needed to cover all vertices of $G$. In this paper, we determine
exact values of isometric path numbers of complete $r$-partite
graphs and Cartesian products of $2$ or $3$ complete graphs.

\bigskip

\bigskip


\end{abstract}

\newpage

\section{Introduction}

An {\it isometric path} between two vertices in a graph $G$ is a
shortest path joining them. The {\it isometric path number} of
$G$, denoted by $\ip(G)$, is the minimum number of isometric paths
required to cover all vertices of $G$. This concept has a close
relationship with the game of cops and robbers described as
follows.

The game is played by two players, the {\it cop} and the
{\it robber}, on a graph. The two players move alternatively,
starting with the cop. Each player's first move consists of
choosing a vertex at which to start. At each subsequent move, a
player may choose either to stay at the same vertex or to move to
an adjacent vertex. The object for the cop is to catch the robber,
and for the robber is to prevent this from happening. Nowakowski
and Winkler \cite{nw83} and Quilliot \cite{q78} independently
proved that the cop wins if and only if the graph can be reduced
to a single vertex by successively removing pitfalls, where a {\it
pitfall} is a vertex whose closed neighborhood is a subset of the
closed neighborhood of another vertex.

As not all graphs are
cop-win graphs, Aigner and Fromme \cite{af84} introduced the
concept of {\it cop-number} of a general graph $G$, denoted by
$c(G)$, which is the minimum number of cops needed to put into the
graph in order to catch the robber . On the way to give an upper
bound for the cop-numbers of planar graphs, they showed that a
single cop moving on an isometric path $P$ guarantee that after a
finite number of moves the robber will be immediately caught if he
moves onto $P$. Observing this fact, Fitzpatrick \cite{f97} then
introduced the concept of isometric path cover and pointed out
that $c(G) \le \ip(G)$.

The isometric path number of the Cartesian product $P_{n_1} \Box
P_{n_2} \Box \ldots \Box P_{n_r}$ has been studied in the
literature. Fitzpatrick \cite{f99} gave bounds for the case when
$n_1=n_2=\ldots=n_r$. Fisher and Fitzpatrick \cite{ff01} gave
exact values for the case $r=2$. Fitzpatrick et al \cite{fnhc01}
gave a lower bound, which is in fact the exact value if $r+1$ is a
power of $2$, for the case when $n_1=n_2=\ldots=n_r=2$. Pan and
Chang \cite{pc} gave a linear-time algorithm to solve the
isometric path problem on block graphs.

In this paper we determine exact values of isometric path numbers
of all complete $r$-partite graphs and Cartesian products of $2$
or $3$ complete graphs. Recall that a {\it complete $r$-partite
graph} is a graph whose vertex set can be partitioned into
disjoint union of $r$ nonempty parts, and two vertices are
adjacent if and only if they are in different parts. We use
$K_{n_1,n_2,\ldots,n_r}$ to denote a complete $r$-partite graph
whose parts are of sizes $n_1, n_2, \ldots, n_r$, respectively. A
{\it Hamming graph} is the Cartesian product of complete graphs,
which is the graph $K_{n_1} \Box K_{n_2} \Box \ldots \Box K_{n_r}$
with vertex set $$
   V(K_{n_1} \Box K_{n_2} \Box \ldots \Box K_{n_r}) =
   \{ (x_1, x_2, \ldots, x_r): 0 \le x_i < n_i  \mbox{ for } 1 \le
   i \le r \}
$$
and edge set $E(K_{n_1} \Box K_{n_2} \Box \ldots \Box K_{n_r})$ is
$$
   \{ (x_1, x_2, \ldots, x_r) (y_1, y_2, \ldots, y_r): x_i=y_i \in
   V(K_i) \mbox{ for all } i \mbox{ except just one } x_j \ne y_j
   \}.
$$

\section{Complete {\boldmath $r$}-partite graphs}

The purpose of this section is to determine exact values of the
isometric path numbers of all complete $r$-partite graphs.

Suppose $G$ is the complete $r$-partite graph
$K_{n_1,n_2,...,n_r}$ of $n$ vertices, where $r \ge 2$, $n_1 \ge
n_2 \ge \ldots \ge n_r$ and $n=n_1+n_2+\ldots+n_r$. Let $G$ has
$\alpha$ parts of odd sizes. We notice that every isometric path
in $G$ has at most $3$ vertices. Consequently,
$$
    \ip(G) \ge \lceil n/3 \rceil.
$$ Also, for any path of $3$ vertices in an isometric path cover
$\cal C$, two end vertices of the path is in a part of $G$ and the
center vertex in another part. In case when two paths of $3$
vertices in $\cal C$ have a common end vertex, we may replace one
by a path of $2$ vertices. And, a path of $1$ vertex can be
replaced by a path of $2$ vertices. So, without loss of
generality, {\it we may only consider isometric path covers in
which every path is of $2$ or $3$ vertices, and two $3$-vertices
paths have different end vertices.}

\begin{lem}             \label{Lemma 1}
If $3 n_1 > 2 n$,
then $\ip(G) = \lceil n_1/2 \rceil$.
\end{lem}
{\bf Proof.}
First,
$\ip(G) \ge \lceil n_1/2 \rceil$ since every
isometric path contains at most two vertices in the first part.

On the other hand, we use an induction on $n-n_1$ to prove that
$\ip(G) \le \lceil n_1/2 \rceil$. When $n-n_1=1$, we have
$G=K_{n-1,1}$. In this case, it is clear that $\ip(G) \le \lceil
n_1/2 \rceil$. Suppose $n-n_1 \ge 2$ and the claim holds for
$n'-n_1' < n-n_1$. Then we remove two vertices from the first part
and one vertex from the second part to form an isometric 3-path
$P$. Since $3n_1 > 2n$, we have $n_1-2 > 2(n-n_1-1) > 0$ and so
$n_1-2 > n_2$. Then, the remaining graph $G'$ has $r' \ge 2$,
$n_1' = n_1-2$ and $n' = n-3$. It then still satisfies $3n_1' >
2n'$. As $n'-n_1' = n-n_1-1$, by the induction hypothesis,
$\ip(G') \le \lceil n_1'/2 \rceil$ and so $\ip(G) \le \lceil
n_1'/2 \rceil + 1 = \lceil n_1/2 \rceil$. \qed

\begin{lem}            \label{Lemma 2}
If $3\alpha > n$, then $\ip(G) = \lceil (n+\alpha)/4 \rceil$.
\end{lem}
{\bf Proof.} Suppose $\cal C$ is an optimum isometric path cover
with $p_2$ paths of 2 vertices and $p_3$ paths of 3 vertices. Then
$$
     2 p_2 + 3 p_3 \ge n.
$$
Notice that there are at most $n-\alpha$ vertices in $G$ can be paired up
as the end vertices of the 3-paths in $\cal P$.
Hence $p_3 \le (n-\alpha)/2$ and so
$$
     2 p_2 + 2 p_3 \ge n - (n-\alpha)/2 = (n+\alpha)/2
     \mbox{ \ or \ }
     \ip(G) = p_2 + p_3 \ge \lceil (n+\alpha)/4 \rceil.
$$

On the other hand, we use an induction on $n-\alpha$ to prove that
$\ip(G) \le \lceil (n+\alpha)/4 \rceil$. When $n-\alpha \le 1$, we
have $n=\alpha$ and $G$ is the complete graph of order $n$. So,
$\ip(G)=\lceil n/2 \rceil = \lceil (n+\alpha)/4 \rceil$. Suppose
$n-\alpha \ge 2$ and the claim holds for $n'-\alpha' < n-\alpha$.
In this case, $3\alpha > n \ge \alpha+2$ which implies $\alpha >
1$ and $n>3$. Then we may remove two vertices from the first part
of and one vertex form an odd part other than the first part to
form a isometric 3-path $P$ of $G$. The remaining graph $G'$ has
$n'=n-3$ and $\alpha'=\alpha-1$. It then satisfies $3\alpha' >
n'$. Notice that $r' \ge 2$ unless $G = K_{2,1,1}$ in which $n=4$
and $\alpha=2$ imply $\ip(G)=2=\lceil (n+\alpha)/4\rceil$. By the
induction hypothesis, $\ip(G') \le \lceil (n'+\alpha')/4 \rceil$
and so $\ip(G) \le \lceil (n'+\alpha')/4 \rceil + 1 = \lceil
(n+\alpha)/4 \rceil$. \qed

\begin{lem}               \label{Lemma 3}
If $3n_1 \le 2n$ and $3\alpha \le n$,
then $\ip(G) = \lceil n/3 \rceil$.
\end{lem}
{\bf Proof.}
Since every isometric path in $G$ has at most $3$ vertices,
$\ip(G) \ge \lceil n/3 \rceil$.

On the other hand, we use an induction on $n$ to prove that
$\ip(G) \le \lceil n/3 \rceil$. When $n\le 8$, by the assumptions
that $3n_1 \le 2n$ and $3\alpha \le n$ we have $G \in \{ K_{2,1}$,
$K_{2,2}$, $K_{3,2}$, $K_{2,2,1}$, $K_{4,2}$, $K_{4,1,1}$,
$K_{3,3}$, $K_{3,2,1}$, $K_{2,2,2}$, $K_{2,2,1,1}$, $K_{4,3}$,
$K_{4,2,1}$, $K_{3,2,2}$, $K_{2,2,2,1}$, $K_{5,3}$, $K_{5,2,1}$,
$K_{4,4}$, $K_{4,3,1}$, $K_{4,2,2}$, $K_{4,2,1,1}$, $K_{3,3,2}$,
$K_{3,2,2,1}$, $K_{2,2,2,2}$, $K_{2,2,2,1,1} \}$. It is
straightforward to check that $\ip(G) \le \lceil n/3 \rceil$.

Suppose $n \ge 9$ and the claim holds for $n' < n$.
We remove two vertices from the first part
and one vertex from the $j$th part to form an isometric $3$-path $P$ for $G$,
where $j$ is the largest index such that $j \ge 2$ and $n_j$ is odd
(when $n_i$ are even for all $i\ge 2$, we choose $j=r$).
Then,
the remaining subgraph $G'$ has $n'=n-3$ and
$\alpha'=\alpha-1$ or $\alpha' \le 2$.
Therefore,
$3\alpha \le n$ and $n \ge 9$ imply that $3\alpha' \le n'$ in any case.
We shall prove that $3n_1' \le 2n'$ according to the following cases.

{\bf Case 1.} $n_1\ge n_2+2$.

In this case, $n_1 - 2 \ge n_2 \ge n_i$ for all $i \ge 2$ and so
$n_1'=n_1-2$. Therefore, $3 n_1' = 3(n_1-2) \le 2(n-3) = 2n'$.

{\bf Case 2.} $n_1 \le n_2+1$ and $n_2 \le 4$.

In this case, $n_1' \le n_2 \le 4$ and $n' \ge 6$.  Then, $3 n_1'
\le 12 \le 2n'$.

{\bf Case 3.} $n_1 \le n_2+1$ and $n_2 \ge 5$ and $r=2$.

In this case, $n_1' \le n_2 -1$ and $n'=n-3 =n_1+n_2-3 \ge
2n_2-3$.  Then, $3n_1' \le 3n_2 - 3 \le 4n_2-8 <2n'$.

{\bf Case 4.} $n_1 \le n_2+1$ and $n_2 \ge 5$ and $r\ge 3$.

In this case, $n_1' \le n_2$ and $n'=n-3 \ge n_1+n_2+1-3 \ge
2n_2-2$. Then, $3n_1' \le 3n_2 \le 4n_2-5 < 2n'$. \qed

According to Lemma \ref{Lemma 1}, \ref{Lemma 2} and \ref{Lemma 3}, we have
the following theorem.

\begin{thm} \label{Theorem 4}
Suppose $G$ is the complete $r$-partite graph
$K_{n_1,n_2,\ldots,n_r}$ of $n$ vertices with $r\ge 2$, $n_1 \ge
n_2 \ge \ldots\ge n_r$ and $n=n_1 + n_2 + \ldots + n_r$.  If there
are exactly $\alpha$ indices $i$ with $n_i$ odd, then $$
\ip(G)=\left\{\begin{array}{ll}
     \lceil n_1/2 \rceil,        & \mbox{if $3n_1 > 2n$}; \\
     \lceil (n+\alpha)/4 \rceil, & \mbox{if $3\alpha > n$}; \\
     \lceil n/3 \rceil,          & \mbox{if $3\alpha \le n$ and $3n_1 \le 2n$.}
              \end{array}\right.
$$
\end{thm}

\bigskip

In the proofs of the lemmas above, the essential points for the
arguments is not the fact that each partite set of the complete
$r$-partite graph is trivial.  If we add some edges into the graph
but still keep that each partite set can be partitioned into
$\lfloor n_i/2 \rfloor$ pairs of two nonadjacent vertices and $n_i
- 2 \lfloor n_i/2 \rfloor$ vertex, then the same result still
holds.

\begin{cor} \label{Corollary 5}
Suppose $G$ is the graph obtained from the complete $r$-partite
graph $K_{n_1,n_2,\ldots,n_r}$ of $n$ vertices by adding edges
such that each $i$-th part can be partitioned into $\lfloor n_i/2
\rfloor$ pairs of two nonadjacent vertices and $n_i - 2 \lfloor
n_i/2 \rfloor$ vertex, where $r\ge 2$, $n_1 \ge n_2 \ge \ldots\ge
n_r$ and $n=n_1+n_2+\ldots+n_r$. If there are exactly $\alpha$
indices $i$ with $n_i$ odd, then $$
\ip(G)=\left\{\begin{array}{ll}
     \lceil n_1/2 \rceil,        & \mbox{if $3n_1 > 2n$}; \\
     \lceil (n+\alpha)/4 \rceil, & \mbox{if $3\alpha > n$}; \\
     \lceil n/3 \rceil,          & \mbox{if $3\alpha \le n$ and $3n_1 \le 2n$.}
              \end{array}\right.
$$
\end{cor}

\section{Hamming graphs}

This section establishes isometric path numbers of Cartesian
products of 2 or 3 complete graphs.

Suppose $G$ is the Hamming graph $K_{n_1} \Box K_{n_2} \Box \ldots
\Box K_{n_r}$ of $n$ vertices, where $n=n_1 n_2 \ldots n_r$ and
$n_i \ge 2$ for $1 \le i \le r$. We notice that every isometric
path in $G$ has at most $r+1$ vertices.  Consequently, $$
     \ip(G) \ge \lceil n/(r+1) \rceil.
$$
Recall that the vertex set of $K_{n_1} \Box K_{n_2} \Box \ldots
\Box K_{n_r}$ is
$$
   V(K_{n_1} \Box K_{n_2} \Box \ldots \Box K_{n_r}) =
   \{ (x_1, x_2, \ldots, x_r): 0 \le x_i < n_i  \mbox{ for } 1 \le
   i \le r \}.
$$
We first consider the case when $r=2$

\begin{thm} \label{Theorem 6}
If $n_1 \ge 2$ and $n_2 \ge 2$, then $\ip(K_{n_1} \Box K_{n_2})=
\lceil n_1 n_2 /3 \rceil$.
\end{thm}
{\bf Proof.} We only need to prove that $\ip(K_{n_1} \Box K_{n_2})
\le \lceil n_1 n_2 /3 \rceil$.  We shall prove this assertion by
induction on $n_1+n_2$. For the case when $n_1+n_2 \le 6$, the
isometric path covers
\begin{eqnarray*}
 {\cal C}_{2,2} &=& \{ (0,0)(0,1),(1,0)(1,1 )\}, \\
 {\cal C}_{2,3} &=& \{ (0,0)(0,1)(1,1), (0,2)(1,2)(1,0) \}, \\
 {\cal C}_{2,4} &=& \{ (0,0)(0,1)(1,1), (0,2)(1,2)(1,0), (0,3)(1,3) \}  \mbox{ and} \\
 {\cal C}_{3,3} &=& \{ (0,0)(2,0)(2,2), (0,1)(0,2)(1,2), (1,0)(1,1)(2,1) \}
 \end{eqnarray*}
for $K_2 \Box K_2$, $K_2 \Box K_3$, $K_2 \Box K_4$ and $K_3 \Box
K_3$ respectively, gives the assertion.

Suppose $n_1+n_2 \ge 7$ and the assertion holds for
$n_1'+n_2'<n_1+n_2$. For the case when all $n_i \le 4$, without
loss of generality we may assume that $n_1 =4$ and $3 \le n_2 \le
4$. As we can partition the vertex set of $K_{n_1} \Box K_{n_2}$
into the vertex sets of two copies of distance invariant induced
subgraphs $K_2 \Box K_{n_2}$,
$$
    \ip(K_{n_1} \Box K_{n_2})
        \le 2 \ip(K_2 \Box K_{n_2})
        \le 2 \lceil 2n_2/3 \rceil
        =   \lceil n_1 n_2/3 \rceil.
$$
For the case when there is at least one $n_i \ge 5$, say $n_1 \ge
5$, again we can partition the vertex set of $K_{n_1} \Box
K_{n_2}$ into the vertex sets of two distance invariant induced
subgraphs $K_3 \Box K_{n_2}$ and $K_{n_1-3} \Box K_{n_2}$. Then,
$$
    \ip(K_{n_1} \Box K_{n_2})
        \le \ip(K_3 \Box K_{n_2}) + \ip(K_{n_1-3} \Box K_{n_2})
        \le \lceil 3n_2/3 \rceil + \lceil (n_1-3)n_2/3 \rceil
        =   \lceil n_1 n_2/3 \rceil.
$$
\qed

\begin{lem} \label{Lemma 7}
If $n_1, n_2$ and $n_3$ are positive even integers, then
$$\ip(K_{n_1} \Box K_{n_2} \Box K_{n_3}) = n_1 n_2 n_3 /4.$$
\end{lem}
{\bf Proof.} We only need to prove that $\ip(K_{n_1} \Box K_{n_2}
\Box K_{n_3}) \le n_1 n_2 n_3 / 4$. First, the isometric path
cover ${\cal C}_{2,2,2}=\{ (0,0,0)(0,0,1)(0,1,1)(1,1,1),
(1,0,1)(1,0,0)(1,1,0)(0,1,0) \}$ for $K_2 \Box K_2 \Box K_2$
proves the assertion for the case when $n_1=n_2=n_3=2$. For the
general case, as the vertex set of $K_{n_1} \Box K_{n_2} \Box
K_{n_3}$ can be partitioned into the vertex sets of $n_1n_2n_3/8$
copies of distance invariant induced subgraphs $K_2 \Box K_2 \Box
K_2$, $$
    \ip(K_{n_1} \Box K_{n_2} \Box K_{n_3})
        \le (n_1 n_2 n_3 /8) \ip(K_2 \Box K_2 \Box K_2)
        \le  n_1 n_2 n_3/4.
$$
\qed

\begin{lem} \label{Lemma 8}
If $n_3 \ge 3$ is odd, then $\ip(K_2\Box K_2\Box K_{n_3})=n_3+1$.
\end{lem}
{\bf Proof.} First, we claim that $\ip(K_2\Box K_2\Box K_{n_3})\ge
n_3+1$. Suppose to the contrary that the graph can be covered by
$n_3$ isometric paths
$$
    P_i: (x_{i1},x_{i2},x_{i3})  (y_{i1},y_{i2},y_{i3})
         (z_{i1},z_{i2},z_{i3})  (w_{i1},w_{i2},w_{i3}),
$$
$i=1,2,\ldots,n_3$. These paths are in fact vertex-disjoint paths
of $4$ vertices, each contains exactly one type-$j$ edge for
$j=1,2,3$, where an edge $(x_1,x_2,x_3)(y_1,y_2,y_3)$ is type-$j$
if $x_j\ne y_j$. For each $P_i$ we then have $x_{i1}=1-w_{i1}$ and
$x_{i2}=1-w_{i2}$, which imply that $x_{i1}+x_{i2}$ has the same
parity with $w_{i1}+w_{i2}$.  We call the path $P_i$ {\it even} or
{\it odd} when $x_{i1}+x_{i2}$ is even or odd, respectively. Also,
as $P_i$ has just one type-$3$ edge, by symmetric, we may assume
either $x_{i3} \ne y_{i3} = z_{i3} = w_{i3}$ or $x_{i3} = y_{i3}
\ne z_{i3} = w_{i3}$, for which we call $P_i$ type 1-3 or type 2-2
respectively.  For a type 2-2 path $P_i$ we may further assume
that $x_{i1} \ne y_{i1} = z_{i1} = w_{i1}$.

For $0\le x_3 < n_3$, the {\it $x_3$-square} is the set $S(x_3) =
\{ (0,0,x_3), (0,1,x_3), (1,0,x_3)$, $(1,1,x_3) \}$.  Notice that
a type 1-3 path $P_i$ contains 1 vertex in $S(x_{i3})$ and 3
vertices in $S(w_{i3})$, while a type 2-2 path $P_i$ contains 2
vertices in $S(x_{i3})$ and 2 vertices in $S(w_{i3})$. We call a
type 1-3 path $P_i$ is {\it adjacent to} another type 1-3 path
$P_j$ if the last 3 vertices of $P_i$ and the first vertex of
$P_j$ form a square.  This defines a digraph $D$ whose vertices
are all type 1-3 paths, in which each vertex has out-degree one
and in-degree at most one.  In fact, each vertex then has
in-degree one. In other words, the ``adjacent to" is a bijection.
Consequently, vertices of all type 1-3 paths together form $p$
squares; and so vertices of all type 2-2 paths form the other
$n_3-p$ squares.

Since $x_{i1} \ne y_{i1} = z_{i1} = w_{i1}$ for a type 2-2 path
$P_i$, the first two vertices of a type 2-2 path together with the
first two vertices of another type 2-2 path form a square.  This
shows that there is an even number of type 2-2 paths. Therefore,
there is an odd number of type 1-3 paths.

On the other hand, in a type 1-3 path $P_i$ we have
$x_{i_1}+x_{i_2}=y_{i_1}+y_{i_2}$ has the different parity with
$z_{i_1}+z_{i_3}$, and the same parity with $w_{i_1}+w_{i_2}$.  So
it is adjacent to a type 1-3 path whose parity is the same as
$z_{i_1}+z_{i_2}$.  That is, a type 1-3 path is adjacent to a type
1-3 path of different parity. Therefore, the digraph $D$ is the
union of some even directed cycle.  This is a contradiction to the
fact that there is an odd number of type 1-3 paths.

The arguments above prove that $\ip(K_2\Box K_2\Box K_{n_3})\ge
n_3+1$. On the other hand, since the vertex set of $K_2\Box
K_2\Box K_{n_3}$ is the union of the vertex sets of $(n_3+1)/2$
copies of $K_2\Box K_2\Box K_2$, by the cover ${\cal C}_{2,2,2}$
in the proof of Lemma \ref{Lemma 7}, we have $\ip(K_2\Box K_2\Box
K_{n_3})\le n_3+1$. \qed

\begin{thm}
If all $n_i\ge 2$, then $\ip(K_{n_1}\Box K_{n_2}\Box
K_{n_3})=\lceil n_1n_2n_3/4 \rceil$ except for the case when two
$n_i$ are $2$ and the third is odd. In the exceptional case,
$\ip(K_{n_1}\Box K_{n_2}\Box K_{n_3})=n_1n_2n_3/4+1$.
\end{thm}
{\bf Proof.} The exceptional case holds according to Lemma
\ref{Lemma 8}.

For the main case, by Lemma \ref{Lemma 7}, we may assume that at
least one $n_i$ is odd. Again, we only need to prove that
$\ip(K_{n_1}\Box K_{n_2}\Box K_{n_3})\le\lceil n_1n_2n_3/4
\rceil$. We shall prove the assertion by induction on
$\sum_{i=1}^3n_i$. For the case when $\sum_{i=1}^3 n_i \le 10$,
the following isometric path covers for $K_2\Box K_3\Box K_3$,
$K_2\Box K_3\Box K_4$, $K_3\Box K_3\Box K_3$ and $K_3\Box K_3\Box
K_4$, respectively, prove the assertion:
\[\begin{array}{lll}{\cal C}_{2,3,3}&=&\{(0,1,1)(0,1,0)(0,0,0)(1,0,0), \ (0,2,2)(0,2,0)(1,2,0)(1,1,0),\\
&& \ \ (0,2,1)(1,2,1)(1,1,1), \ (0,0,2)(0,1,2)(1,1,2),\\ && \ \
(0,0,1)(1,0,1)(1,0,2)(1,2,2)\};
\end{array}\]
\[\left( \begin{array}{lll} {\rm Let } \ {\cal C}_{2,3,3}^*&=&{\cal
C}_{2,3,3}\backslash\{(0,2,1)(1,2,1)(1,1,1),(0,0,2)(0,1,2)(1,1,2)\} \cup\\
&& \{(0,2,1)(1,2,1)(1,1,1)(1,1,3),(0,0,2)(0,1,2)(1,1,2)(1,1,4)\}.
\end{array} \right)\]
\[\begin{array}{lll}{\cal C}_{2,3,4}&=&\{(0,1,1)(0,1,0)(0,0,0)(1,0,0), \ (0,2,1)(0,2,0)(1,2,0)(1,1,0),\\
&& \ \ (0,2,3)(0,2,2)(1,2,2)(1,1,2), \
(0,1,3)(0,1,2)(0,0,2)(1,0,2),\\ && \ \
(0,0,1)(1,0,1)(1,1,1)(1,1,3), \ (1,2,1)(1,2,3)(1,0,3)(0,0,3)\};
\end{array}\]
\[\begin{array}{lll}{\cal C}_{2,3,5}&=&{\cal C}_{2,3,3}^*\cup\{(0,1,4)(0,1,3)(0,2,3)(1,2,3), \ (0,0,3)(0,0,4)(0,2,4)(1,2,4),\\
&& \hskip 1.5cm (1,0,3)(1,0,4)\};
\end{array}\]
\[\begin{array}{lll}{\cal C}_{3,3,3}&=&\{(0,0,0)(0,2,0)(1,2,0)(1,2,1), \ (1,1,0)(2,1,0)(2,2,0)(2,2,1),\\
&& \ \ (0,2,1)(0,1,1)(1,1,1)(1,1,2), \ (1,0,1)(2,0,1)(2,1,1)(2,1,2),\\
&& \ \ (0,1,0)(0,1,2)(0,2,2)(1,2,2), \ (0,0,1)(0,0,2)(2,0,2)(2,2,2),\\
&& \ \ (1,0,2)(1,0,0)(2,0,0)\};
\end{array}\]
\[\begin{array}{lll}{\cal C}_{3,3,4}&=&\{(0,0,0)(0,2,0)(1,2,0)(1,2,1), \ (1,1,0)(2,1,0)(2,2,0)(2,2,1),\\
&& \ \ (0,2,1)(0,1,1)(1,1,1)(1,1,2), \ (1,0,1)(2,0,1)(2,1,1)(2,1,2),\\
&& \ \ (0,1,0)(0,1,2)(0,2,2)(1,2,2), \ (0,0,2)(2,0,2)(2,2,2)(2,2,3),\\
&& \ \ (0,1,3)(1,1,3)(1,0,3)(1,0,2), \ (1,0,0)(2,0,0)(2,0,3)(2,1,3),\\
&& \ \ (0,0,1)(0,0,3)(0,2,3)(1,2,3)\}.
\end{array}\]
Suppose $\sum_{i=1}^3 n_i \ge 11$ and the assertion holds for
$\sum_{i=1}^3n_i'<\sum_{i=1}^3n_i$. We shall consider the
following cases.

For the case when there is some $i$, say $i=3$, such that $n_3 \ge
7$ or $n_3=6$ with all $n_j \ge 3$, we have $
   \ip(K_{n_1} \Box K_{n_2} \Box K_{n_3}) \le
   \ip(K_{n_1} \Box K_{n_2} \Box K_4) + \ip(K_{n_1} \Box K_{n_2} \Box
     K_{n_3-4}) \le \lceil n_1n_24/4 \rceil + \lceil
     n_1n_2(n_3-4)/4 \rceil = \lceil n_1n_2n_3/4 \rceil.
$

For the case when some $n_i$, say $n_3$, is equal to $4$, we may
assume $n_1\ge n_2$ and so $n_1 \ge 4$.  Then $
 \ip(K_{n_1} \Box K_{n_2} \Box K_4) \le \ip(K_2 \Box
  K_{n_2} \Box K_4) + \ip(K_{n_1-2} \Box K_{n_2} \Box
  K_4) = \lceil 2 n_2 4/4 \rceil + \lceil (n_1-2)n_2 4/4 \rceil
  = \lceil n_1n_2n_3/4 \rceil.
$

There are $6$ remaining cases. The following isometric path covers
prove the assertion for $K_2\Box K_3\Box K_6$, $K_2\Box K_5\Box
K_5$ and $K_3\Box K_5\Box K_5$, respectively:
\[\begin{array}{lll}{\cal C}_{2,3,6}&=&{\cal C}_{2,3,3}^*\cup\{(0,0,4)(0,0,3)(1,0,3)(1,2,3), \ (0,1,3)(0,1,4)(0,2,4)(1,2,4),\\
&& \hskip 1.5cm (0,2,3)(0,2,5)(1,2,5)(1,1,5), \
(0,1,5)(0,0,5)(1,0,5)(1,0,4)\};
\end{array}\]
\[\begin{array}{lll}{\cal C}_{2,5,5}&=&{\cal C}_{2,3,5}\backslash\{(1,0,3)(1,0,4)\}\cup\\
&&  \{(0,4,1)(0,4,0)(0,3,0)(1,3,0), \
(1,4,0)(1,4,1)(1,3,1)(0,3,1),\\ && \ \
(0,4,3)(0,4,2)(0,3,2)(1,3,2), \ (1,4,2)(1,4,3)(1,3,3)(0,3,3),\\ &&
\ \ (1,0,3)(1,0,4)(1,4,4), \ (0,4,4)(0,3,4)(1,3,4)\};
\end{array}\]
\[\begin{array}{lll}{\cal C}_{3,5,5}&=&{\cal C}_{2,3,5}\backslash\{(1,0,3)(1,0,4)\}\cup\\
&& \{(0,4,0)(2,4,0)(2,0,0)(2,0,1), \
(0,3,0)(2,3,0)(2,1,0)(2,1,1),\\ && \ \
(0,4,1)(0,3,1)(1,3,1)(1,3,0), \ (1,4,0)(1,4,1)(2,4,1)(2,2,1),\\ &&
\ \ (1,0,3)(2,0,3)(2,2,3)(2,2,0), \
(1,0,4)(2,0,4)(2,3,4)(2,3,1),\\ && \ \
(0,3,2)(2,3,2)(2,1,2)(2,1,3), \ (0,4,4)(0,4,2)(2,4,2)(2,0,2),\\ &&
\ \ (0,4,3)(1,4,3)(1,3,3)(1,3,2), \
(0,3,3)(2,3,3)(2,4,3)(2,4,4),\\ && \ \
(0,3,4)(1,3,4)(1,4,4)(1,4,2), \ (2,2,2)(2,2,4)(2,1,4)\}.
\end{array}\]
The other $3$ cases follows from the following inequalities:
$$
   \ip(K_2 \Box K_5 \Box K_6) \le \ip(K_2 \Box K_3 \Box K_6) + \ip(K_2 \Box
    K_2 \Box K_6) \le 9+6 = 15,
$$
$$
   \ip(K_3 \Box K_3 \Box K_5) \le \ip(K_3 \Box K_3 \Box K_2) + \ip(K_3 \Box K_3 \Box K_3)
   \le 5+7 = 12,
$$
$$
   \ip(K_5 \Box K_5 \Box K_5) \le \ip(K_5 \Box K_5 \Box K_3) +
   \ip(K_5 \Box K_5 \Box K_2) \le 19+13 =32.
$$
\qed

\end{document}